\documentclass{amsart}
\usepackage{amssymb}  
\usepackage{amsthm} 
\usepackage{amsmath}
\usepackage{amsfonts}
\usepackage{graphicx}
\usepackage[utf8]{inputenc}
\usepackage[all]{xy}
\usepackage[linktocpage=true,colorlinks=true, allcolors=blue]{hyperref}
\author{José Juan-Zacarías}
\address{Instituto de Matemáticas Unidad
 Cuernavaca, Av. Universidad s/n. Col. Lomas de Chamilpa Código
 Postal 62210, Cuernavaca, Morelos.}
\email{jose.juan@im.unam.mx}

\date{October 19, 2018}
\begin{document}
\title{Baker-Gross theorem revisited}
\maketitle

\begin{abstract} 
F. Gross conjectured that any meromorphic solution of the Fermat
Cubic $F_3\colon\ x^3+y^3=1$ are elliptic functions composed with
entire functions. The conjecture was solved affirmatively first by I.
N. Baker who found explicit formulas of those elliptic functions and
later F. Gross gave another proof proving that in fact one of them
uniformize the Fermat cubic. In this paper we give an alternative
proof of the Baker and Gross theorems. With our method we obtain
other analogous formulas. Some remarks on Fermat curves of higher
degree is given.
\end{abstract}

\section*{Introduction}
Consider the Fermat cubic
\begin{equation}\label{Fermat-cubic}
  F_3\colon \ x^3+y^3=1.
\end{equation}

This algebraic curve defines an elliptic curve, i.e., a compact
Riemann surface of genus 1 (taking the zeros in $\mathbb{CP}^2$ of
its homogenization). A meromorphic solution of this equation is, by
definition, a pair of meromorphic functions in the plane such that
$f^3+g^3=1$. In his paper \cite{Gross} F. Gross conjectures that any
meromorphic solution of the Fermat cubic is obtained by composing
elliptic functions with entire functions. The conjecture was solved
affirmatively by I. N. Baker in \cite{Baker}. He proved that any
solution is the composition of  the following elliptic  functions
with an entire function:
 
\begin{equation}\label{elliptic-solutions}
f(z)=\frac{1}{2\wp(z)}\left(1-3^{-1/2}\wp'(z)\right), \quad
g(z)=\frac{1}{2\wp(z)}\left(1+3^{-1/2}\wp'(z)\right), \end{equation}
where $\wp$ is the Weierstrass elliptic function satisfying
$(\wp')^2=4\wp^3-1$. In what follows we denote by $\Lambda'$ the
lattice  in $\mathbb C$ that defines this $\wp$. In particular these
functions are solutions of the Fermat cubic but these formulas differ
from the analogous that appear in \cite{Gross}, \cite{Gross-erratum},
which seem to contain an error. Later, F. Gross gave another proof in
\cite{Gross-II}, proving in fact that the function $f$ in
(\ref{elliptic-solutions}) gives a uniformization of the Fermat cubic
(\ref{Fermat-cubic}). In our context we formulate the previous
results in the following theorem:\\ \par \noindent \textbf{Theorem
(Baker-Gross).} \textit{Let $\Lambda'$
  and $\wp$ be as above. Then the map
$\mathbb{C}/\Lambda'\to F_3$ given in affine coordinates by
\begin{equation}\label{map-uniform}
z\mapsto \left(\frac{1}{2\wp(z)}\left(1-3^{-1/2}\wp'(z)\right),
\frac{1}{2\wp(z)}\left(1+3^{-1/2}\wp'(z)\right)\right)
\end{equation}
is a biholomorphism between the two elliptic curves. Then by the
lifting property of coverings, any pair of functions $F$ and $G$,
which are meromorphic in the plane and satisfy (\ref{Fermat-cubic})
have the form:
\begin{equation}\label{general-form}
  F=\frac{1}{2\wp(\alpha)}\left(1-3^{-1/2}\wp'(\alpha)\right),
  \quad G=\frac{1}{2\wp(\alpha)}\left(1+3^{-1/2}\wp'(\alpha)\right),
\end{equation}
where $\alpha$ is an entire function.}\\

\par In this paper we give a proof of this theorem by using Riemann
surface theory and by using an explicit map from a Weierstrass normal
form to the Fermat cubic. Our proof could clarify the nature of the
previous formulas, which are not obvious. Also, by this method, other
formulas analogous to (\ref{map-uniform}) and (\ref{general-form})
are obtained (see (\ref{map-uniform2}) and
(\ref{general-form2})).

\par In Section 1 we recall some basic facts about elliptic curves
and compute a Weierstrass normal form of the Fermat cubic, and the
corresponding isomorphism as well.  In the next section we prove the
main theorem. Finally, in the last section we give some remarks on
Fermat curves of higher degree.

\par Recently, N. Steinmetz communicated to the author another proof
of the Gross conjecture in \cite{Steinmetz} (\S 2.3.5 pp. 56-57) by
using Nevanlinna theory. He proved without reference to the
Uniformization Theorem the following:\\

\par \noindent \textbf{Theorem (Steinmetz).} \textit{Suppose that
non-constant meromorphic functions $f$ and $g$ parametrize the
algebraic curve
\begin{equation*}
F\colon\ x^n+y^m=1\quad (n\geq m\geq 2)
\end{equation*}
with $\frac{1}{m}+\frac{1}{n}<1$. Then $(m,n)$ equals $(4,2)$ or
$(3,3)$ or $(3,2)$. In any case $f$ and $g$ are given by
\begin{equation*}
 f=E\circ \psi\quad \text{and}\quad g=\sqrt[m-1]{E'}\circ \psi,
\end{equation*}
where $E$ is an elliptic function satisfying
\begin{equation*}
 E'^2=1-E^4,\quad E'^3=(1-E^3)^2\quad and \quad E'^2=1-E^3,
\end{equation*}
respectively, and $\psi$ is any non-constant entire function.}\\

The present paper contains part of the Undergraduate Thesis of the
author written under the supervision of Dr. Alberto Verjovsky at the
Cuernavaca Branch of the Institute of Mathematics of the National
Autonomous University of Mexico (UNAM).
\section{The normal form of the Fermat cubic}

\subsection{Basic facts on elliptic curves} A complex elliptic curve
$X$ is by definition a compact Riemann surface of genus 1. The
Plücker formula tells us that a non-singular projective curve of
degree 3 in $\mathbb{CP}^2$ is a Riemann surface of genus 1 
i.e., an elliptic curve. The reciprocal is also true and we
will briefly discuss it. For this, we recall the uniformization
theorem and the Weierstrass normal form.

\par The \emph{Uniformization Theorem} says that every simply
connected Riemann surface is conformally equivalent to one of the
three Riemann surfaces: the Riemann sphere $\overline{\mathbb{C}}$,
the complex plane $\mathbb{C}$, or the open unit disk $\Delta$. This
theorem combined with the theory of covering spaces give us a
classification of Riemann surfaces: every Riemann surface $X$ is
conformally equivalent to a quotient $\tilde{X}/G$, where $\tilde{X}$
is the universal holomorphic cover of $X$ (hence isomorphic to one of
the three previous Riemann surfaces) and $G$ is a subgroup of
holomorphic automorphisms of $\tilde{X}$ which acts on $\tilde{X}$
free and properly discontinuously. In particular, when the Riemann
surface is of genus 1, it has the complex plane as its universal
holomorphic cover, then $X$ is conformally equivalent to
$\mathbb{C}/\Lambda$, for some lattice $\Lambda\subset \mathbb{C}$.
For an introduction to Riemann surfaces and a proof of the
uniformization theorem see \cite{Forster}.

\par The homogeneous polynomial with complex coefficients

\begin{equation}\label{homog-cubic}
	Y^2Z-4X^3+g_2XZ^2+g_3Z^3,
\end{equation}
obtained by homogenization of the polynomial
\begin{equation}\label{no-hom-cubic}
y^2=4x^3-g_2x-g_3,
\end{equation}
defines an non-singular curve if and only if the discriminant
$\Delta=g_2^3-27g_3^2$ does not vanish. Hence, (\ref{homog-cubic})
defines an elliptic curve if and only if $\Delta\neq 0$. We call a
\emph{Weierstrass normal form} of an elliptic curve $X$ an elliptic
curve given by an equation of the form (\ref{homog-cubic}) which is
isomorphic as a Riemann surface to $X$.

\par Recall also that given a lattice $\Lambda \subset \mathbb{C}$ we
can associate the Weierstrass elliptic function $\wp$ or
$\wp_\Lambda$ given by the series:
\begin{equation}
  \wp(z)=\frac{1}{z^2}+\sum_{\omega\in
    \Lambda^*}\left(\frac{1}{(z+\omega)^2}-\frac{1}{\omega^2}\right).
\end{equation}

This function satisfies the differential equation \begin{equation}
(\wp')^2=4\wp^3-g_2\wp-g_3, \end{equation} where $g_2$ and $g_3$ are
constants depending on $\Lambda$ given by:
\begin{equation*}
 g_2= 60 \sum_{\omega\in \Lambda^*}\frac{1}{\omega^4}, \quad
 g_3=140 \sum_{\omega\in \Lambda^*}\frac{1}{\omega^6},
\end{equation*}
satisfying $\Delta=g_2^3-27g_3^2\neq 0$. Thus this function
gives us a map $\Psi\colon \mathbb{C}/\Lambda\to E$, in affine
coordinates given by: \begin{equation}\label{Psi}
\Psi(z)=(\wp(z),\wp'(z)), \end{equation} from $\mathbb{C}/\Lambda$ to
the elliptic curve $E\colon y^2=4x^3-g_2x-g_3$. This map is an
biholomorphism which sends $\Lambda$ to the point at infinity
$[0:1:0]$.

\par From the previous results and the Uniformization Theorem we can
conclude that every elliptic curve has a Weierstrass normal form.
\par Also, it is true that  given a non-singular equation
(\ref{no-hom-cubic}), there exists a lattice $\Lambda$ with the same
constants $g_2$ and $g_3$. For more information, refer to \cite[p.
176]{Silverman}.

\subsection{Computing the Weierstrass normal form of the Fermat
cubic}\label{computing} Although a Weierstrass normal form is in
general difficult to compute starting from an abstract Riemann
surface of genus 1, the case of the Fermat cubic is relatively easy
by choosing suitable changes of variables. Since this process will be
applied to other Fermat curves in Section \ref{Remarks}, we
describe it step-by-step below:

\begin{enumerate}
\item[1.]  Change $(x,y)$ to $(x-y,x+y)$ in order to eliminate the
cubic term $y^3$. Obtaining:
\begin{equation*}
  E_1\colon\ 2x^3+6xy^2=1.
\end{equation*}
\item[2.] Change $(x,y)$ to $(1/x,y/x)$  to get: 
\begin{equation*}
  E_2\colon\ 2+6y^2=x^3.
\end{equation*}
\item[3.] At this point, we could use any change of variables for
which the coefficient of $y^2$ is 1 and the coefficient of $x^3$ is
$4$, for instance with $(x, y/\sqrt{24})$ we obtain the case $g_2=0$
and $g_3=8$:
\begin{equation*}
E_3\colon\ y^2=4x^3-8.
\end{equation*}
\end{enumerate}

Observe that we obtain a map from the curve obtained in the change of 
variable to the original curve. For example in step 1 we obtain $E_1\to 
F_3$, $(x,y)\mapsto (x-y,x+y)$. Then, the maps associated to the 
previous changes of variables are:
\begin{eqnarray}\label{Phi-}
E_3\to E_2 \quad \quad & E_2\to E_1 & \quad E_1\to F_{3}\\
\nonumber \quad (x,y)\mapsto \left(x,\frac{y}{\sqrt{24}}\right), &\quad
(x,y) \mapsto \left(\frac{1}{x},\frac{y}{x}\right), &\quad (x,y)\mapsto
(x-y,x+y).
\end{eqnarray}

The inverse maps are (in the reverse order, respectively):
\begin{eqnarray}\label{Phi}
 F_3\to E_1 \quad\quad&\quad E_1\to E_2 &\quad E_2\to E_3\\
\nonumber (x,y)\mapsto \left(\frac{y+x}{2},\frac{y-x}{2}\right), &\quad
(x,y)\mapsto \left(\frac{1}{x},\frac{y}{x}\right), &\quad
(x,y)\to (x,\sqrt{24}y).
\end{eqnarray}

So in each step we have a birrational isomorphism between these
non-singular algebraic curves, hence a biholomorphism between their
Riemann surfaces. So we obtain, composing the maps of (\ref{Phi}) and
(\ref{Phi-}), respectively, the biholomorphisms $\Phi\colon F_3\to
E_3$ and $\Phi^{-1}\colon E_3\to F_3$:
\begin{eqnarray}\label{Phis}
 \Phi(x,y) & = & \left(\frac{2}{y+x},\sqrt{24}\frac{y-x}{y+x}\right),\\
 \nonumber \Phi^{-1}(x,y) & = &\left(\frac{1}{x}-\frac{y}{\sqrt{24}x},
 \frac{1}{x}+\frac{y}{\sqrt{24}x}\right).
\end{eqnarray}

\section{Proof of the Baker-Gross theorem}
From the previous explicit formulas the Baker-Gross theorem follows
easily. Consider $\Lambda$ associated to $g_2=0$ and $g_3=8$ and
consider the biholomorphism $\Psi:\mathbb{C}/\Lambda\to E_3$ defined
in (\ref{Psi}), then the composition $\Phi^{-1}\circ\Psi\colon
\mathbb{C}/\Lambda\to F_3$ is a biholomorphism,
\begin{equation}\label{map-uniform2}
  \Phi^{-1}\circ\Psi (z)=
  \left(\frac{1}{\wp(z)}-\frac{1}{\sqrt{24}}\frac{\wp'(z)}{\wp(z)},
    \frac{1}{\wp(z)}+\frac{1}{\sqrt{24}}\frac{\wp'(z)}{\wp(z)}\right).
\end{equation}
where $\wp$ satisfies $(\wp')^2=4\wp^3-8$.
\par If we continue from step 3 applying the change of variables
$(2x,\sqrt{2^3}y)$ we obtain the curve $E_3'\colon y^2=4x^3-1$ and the
map $\overline{\Phi}=\Phi^{-1}(2x,\sqrt{2^3}y):E_3'\to F_3$
\begin{eqnarray}
\overline{\Phi}(x,y) & = & \Phi^{-1}(2x,\sqrt{2^3}y)\\
\nonumber & =
&\left(\frac{1}{2x}-\frac{\sqrt{2^3}y}{2\sqrt{24}x},\frac{1}{2x}+
\frac{\sqrt{2^3}y}{2\sqrt{24}x}\right)\\
\nonumber & = &
\left(\frac{1}{2x}\left(1-\frac{y}{\sqrt{3}x}\right),\frac{1}{2x}
\left(1+\frac{y}{\sqrt{3}x}\right) \right),
\end{eqnarray}
and taking $\Lambda'$ associated to $g_2=0$ and $g_3=1$, and
$\Psi':\mathbb{C}/\Lambda'\to E_3'$ as (\ref{Psi}), composing  this
two isomorphism we obtain the biholomorphism expected in
(\ref{map-uniform}) $\overline{\Phi}\circ \Psi'\colon
\mathbb{C}/\Lambda'\to F_3$:
\begin{equation*}
\overline{\Phi}\circ
\Psi'(z)=\left(\frac{1}{2\wp(z)}\left(1-3^{-1/2}\wp'(z)\right), 
\frac{1}{2\wp(z)}\left(1+3^{-1/2}\wp'(z)\right)\right),
\end{equation*}
where the Weierstrass elliptic  function $\wp$ satisfies here
$(\wp')^2=4\wp^3-1$.

\par On the other hand, let $\pi\colon\mathbb{C}\to
\mathbb{C}/\Lambda'$ be the natural projection, this map is an
unbranched holomorphic covering, then the map
$\overline{\Phi}\circ\Psi'\circ \pi\colon \mathbb{C}\to F_3$ is an
unbranched holomorphic covering as well. Hence, given $F$ and $G$ a
meromorphic solution of the Fermat cubic, the map
$\phi(z)=(F(z),G(z))$ defines a holomorphic map $\phi\colon
\mathbb{C}\to F_3$. Since $\mathbb{C}$ is simply connected $\phi$ has
an holomorphic lifting $\alpha\colon \mathbb{C}\to \mathbb{C}$ with
respect to this covering, i.e., the following diagram commutes:

\begin{equation}
\xymatrix{ & \mathbb{C}\ar[d]^{\overline{\Phi}\circ\Psi'\circ \pi}\\
		\mathbb{C}\ar^{\alpha}[ru]\ar^{\phi}[r] & F_3 }
\end{equation}

Composing with $\alpha$ we obtain
\begin{equation}
  F=\frac{1}{2\wp(\alpha)}\left(1-3^{-1/2}\wp'(\alpha)\right),
  \quad G=\frac{1}{2\wp(\alpha)}\left(1+3^{-1/2}\wp'(\alpha)\right),
\end{equation}
which are the desired formulas. This proves the theorem. \par Note that
we could use the map $\Phi^{-1}\circ \Psi\colon \mathbb{C}/\Lambda\to
F_3$ given in (\ref{map-uniform2}) instead of $\overline{\Phi}\circ
\Psi'$ in the above argument to obtain that any meromorphic solution
of the Fermat cubic is of the form
\begin{equation}\label{general-form2}
F=\frac{1}{\wp(\alpha)}\left(1-\frac{1}{\sqrt{24}}\wp'(\alpha)\right),
\quad
G=\frac{1}{\wp(\alpha)}\left(1+\frac{1}{\sqrt{24}}\wp'(\alpha)\right),  
\end{equation}
where in this case $\wp$ satisfies $(\wp')^2=4\wp^3-8$. We could
obtain similar solutions depending on which factor we choose in step
3, but we can always obtain one from the other by this process.

\section{Some remarks for Fermat curves of higher
degree.}\label{Remarks}

We finalize discussing about the application of the changes of
variables described in \ref{computing} to the Fermat curves of higher
degrees (see (\ref{Fermat-curve})). When the curve is of odd degree
the process give us directly an interesting equation, but when the
degree is even we need to apply a slight modification in step 1.
From these equations we give a meromorphic function on the Fermat
curves.

\subsection{The odd case} The changes of variables in steps 1 and 2
described in \ref{computing} can be applied to any Fermat curve,

\begin{equation}\label{Fermat-curve} F_n\colon x^n+y^n=1,
\end{equation} but in the case of $n$ odd we get an interesting
formula. By a straightforward calculation, following steps 1 and 2,
we find the curve $E_2$:

\begin{equation}\label{new-Fermat}
 E_2\colon\ 2+2\sum_{k=1}^{\frac{n-1}{2}}\binom{n}{2k}y^{2k}=x^n.
\end{equation}

As we did not modify the above steps we get the same
correspondence $\Phi\colon F_n\to E_2$ as in (\ref{Phis}) but
without step 3, so we get in this case:
\begin{eqnarray}\label{Phis2}
 \Phi(x,y)=\left(\frac{2}{y+x},\frac{y-x}{y+x}\right),\\
 \nonumber \Phi^{-1}(x,y)=\left(\frac{1}{x}-\frac{y}{x},
 \frac{1}{x}+\frac{y}{x}\right).
\end{eqnarray}

Note that $E_2$ has an holomorphic involution $I(x,y)=(x,-y)$. It is
easy to check that it is conjugate by $\Phi$ to the canonical
involution of $F_n$, $\overline{I}(x,y)=(y,x)$, i.e., the
following diagram commutes
\begin{equation}
\xymatrix{F_3\ar_\Phi[d]\ar^{\overline{I}}[r] & F_3\ar^{\Phi}[d]\\
E_2\ar^{I}[r] & E_2}
\end{equation}

Note that the projection in the first coordinate  is a meromorphic
function of degree $n-1$ on $E_2$, so composing with $\Phi$  we
obtain the meromorphic function $2/(y+x)$ on $F_n$ of degree $n-1$,
for example in the case $n=3$ we obtain a degree 2 meromorphic
function on the elliptic curve $F_3$.

\subsection{The even case}
Similar formulas can be obtained in the even case by using the change
$(x+\omega y,x+y)$ instead of $(x-y,x+y)$ in the first step, where
$\omega$ is a root of $x^n=-1$, maintaining the other steps without
changes as before. In this case we have
\begin{equation}
E_2\colon 2+\sum_{k=1}^{n-1}\binom{n}{k}(1+\omega^k)y^{k}=x^n,
\end{equation}
and $\Phi:F_n\to E_2$ become
\begin{eqnarray}
\Phi(x,y) & = & \left(\frac{\omega-1}{\omega y-x},\frac{x-y}{\omega
y-x}\right),\\
\nonumber \Phi^{-1}(x,y) & = & \left(\frac{1}{x}+\omega\frac{y}{x},
\frac{1}{x}+\frac{y}{x}\right).
\end{eqnarray}

Similarly as above, the map $(\omega-1)/(\omega y-x)$ is an
meromorphic map of degree $n-1$ on the Fermat curve $F_n$, for $n$
even.

\section*{Acknowledgments} \par I would like to thank my advisor 
Alberto Verjovsky for his constant support and for his encouragement in 
writing this paper. This work was partially supported by PAPIIT 
IN100811.

\end{document}